\documentclass[11pt]{amsart}
\usepackage[sorted-cites]{amsrefs}
\usepackage{verbatim}
\usepackage{color}
\usepackage{amscd}
\renewcommand{\bar}{\overline}

\newcommand{\ph}{\varphi}

\newcommand{\rr}{\mathbb{R}}
\newfont{\fnt}{cmr10 scaled 550}

\newtheorem{theorem}{Theorem}

\newtheorem{lemma}{Lemma}
\newtheorem{cor}{Corollary}
\newtheorem{prop}{Proposition}

\newtheorem{definition}{Definition}

\theoremstyle{remark}
\newtheorem{remark}{Remark}

\font\strange=msbm10

\renewcommand{\epsilon}{\varepsilon}

\renewcommand{\Sigma}{\varSigma}
\newcommand{\R}{{{\mathchoice  {\hbox{$\textstyle{\text{\strange R}}$}}
{\hbox{$\textstyle{\text{\strange R}}$}}
{\hbox{$\scriptstyle  N\kern-0.3em  R$}}
{\hbox{$\scriptscriptstyle  R\kern-0.2em  R$}}}}}

\newcommand{\Z}{{{\mathchoice  {\hbox{$\textstyle{\text{\strange Z}}$}}
{\hbox{$\textstyle{\text{\strange Z}}$}}
{\hbox{$\scriptstyle  Z\kern-0.3em  Z$}}
{\hbox{$\scriptscriptstyle  Z\kern-0.2em  Z$}}}}}

\newcommand{\N}{{{\mathchoice  {\hbox{$\textstyle{\text{\strange N}}$}}
{\hbox{$\textstyle{\text{\strange N}}$}}
{\hbox{$\scriptstyle  N\kern-0.3em  N$}}
{\hbox{$\scriptscriptstyle  N\kern-0.2em  N$}}}}}

\renewcommand{\phi}{\varphi}
\usepackage{amsmath,amsthm,amscd}

\newcommand{\Ric}{\textrm{Ric}}

\begin{document}

\title{Stability and  compactness for complete $f$-minimal surfaces}



\thanks{The first and third authors are  partially supported by CNPq and Faperj of Brazil. The second author is supported by CNPq of Brazil.}

\author{Xu Cheng}
\address{Instituto de Matematica e Estat\'\i stica, Universidade Federal Fluminense,
Niter\'oi, RJ 24020, Brazil, email: xcheng@impa.br}
\author{Tito Mejia}
\address{Instituto de Matematica e Estat\'\i stica, Universidade Federal Fluminense,
Niter\'oi, RJ 24020, Brazil, email: tmejia.uff@gmail.com}
\author[Detang Zhou]{Detang Zhou}
\address{Instituto de Matematica e Estat\'\i stica, Universidade Federal Fluminense,
Niter\'oi, RJ 24020, Brazil, email: zhou@impa.br}


\newcommand{\M}{\mathcal M}

\begin{abstract} Let $(M,\bar{g}, e^{-f}d\mu)$ be a complete metric measure space with Bakry-\'Emery Ricci curvature bounded below by a positive constant.
We  prove that,  in $M$, there is no complete  two-sided   $L_f$-stable  immersed  $f$-minimal hypersurface with finite weighted volume. Further, if $M$ is a $3$-manifold, we  prove a smooth compactness theorem for the space of complete  embedded $f$-minimal surfaces in $M$ with the uniform upper bounds of genus and  weighted volume, which generalizes the compactness theorem for complete self-shrinkers in $\mathbb{R}^3$ by Colding-Minicozzi. 
\end{abstract}

\maketitle
\baselineskip=18pt

\section{Introduction}\label{introduction}

Recall that a self-shrinker  (for mean curvature flow in $\mathbb{R}^{n+1}$) is a hypersurface $\Sigma$ immersed in the Euclidean space $(\rr^{n+1} , g_{can})$ satisfying that $$H=\frac12 \langle x,\nu\rangle,$$
where $x$ is the position vector in $\mathbb{R}^{n+1}$, $\nu$ is the unit normal at $x$,  and $H$ is the mean curvature curvature of  $\Sigma$ at $x$.
Self-shrinkers play an important role in the study of singularity of mean curvature flow and have been studied by many people in recent years. We refer   to  \cite{CM1},  \cite{CM3} and the references therein. In particular, Colding-Minicozzi \cite{CM1} proved  the following compactness theorem for self-shrinkers in $\rr^3$.
\begin{theorem} \cite{CM1} \label{cm-compact} Given an integer $g\ge 0$ and a constant $V>0$, the space $\mathcal{S}(g,V)$ of smooth complete embedded self-shrinkers $\Sigma\subset \rr^3$ with
\begin{itemize}
\item genus at most $g$,
\item $\partial \Sigma=\emptyset$,
  \item Area$(B_R(x_0)\cap\Sigma)\le VR^2 $ for all $x_0\in \rr^3$ and $R>0$
\end{itemize}
is compact.

Namely, any sequence of these has a subsequence that converges in the topology of $C^m$
convergence on compact subsets for any $m \geq 2$.
\end{theorem}
In this paper,  we extend  Theorem \ref{cm-compact} to the space of complete embedded $f$-minimal surfaces  in a $3$-manifold.   
A  hypersurface  $\Sigma$  immersed in a Riemannian manifold $(M,\bar{g})$ is called a $f$-minimal hypersurface if its   mean curvature $H$ satisfies that, for any $p\in \Sigma$,  $$H=\langle \overline{\nabla} f,\nu\rangle,$$  where $f$ is a smooth function defined on $M$, and $ \overline{\nabla} f$ denotes the gradient of $f$ on $M$.  Here are some examples of $f$-minimal hypersurfaces:
\begin{itemize}
\item $f\equiv C$, an $f$-minimal hypersurface is just a minimal hypersurface;
\item self-shrinker $\Sigma$ in $\mathbb{R}^{n+1}$.   $f=\frac{|x|^2}{4}$;
\item Let  $(M, \bar{g}, f)$ be a shrinking gradient Ricci solitons,  i.e.  after a normalization,  $(M,\bar{g},f)$ satisfies the equation  $
\overline{\textrm{Ric}}+\overline{\nabla}^2f=\frac12 \bar{g}$ or equivalently  the Bakry-\'Emery Ricci curvature $\overline{\textrm{Ric}}_f:=\overline{\textrm{Ric}}+\overline{\nabla}^2f= \frac12$.   We may consider $f$-minimal hypersurfaces in $(M, \bar{g}, f)$. In particular,  the previous example: a self-shrinker $\Sigma$ in $\mathbb{R}^{n+1}$ is $f$-minimal in Gauss shrinking soliton $(\mathbb{R}^{n+1}, g_{can}, \frac{|x|^2}{4})$;
\item  $M=\mathbb{H}^{n+1}(-1)$, the hyperbolic space.  Let  $r$ denote the distance function from a fixed point $p\in M$ and  $f(x)=nar^2(x)$, where $a>0$ is a constant. Now $\overline{\textrm{Ric}_f}\ge n(2a-1)$.  The geodesic sphere of radius $r$ centered at $p$ is an $f$-minimal hypersurface if the radius $r$ satisfies $2ar=\coth r$. 

\end{itemize}

An $f$-minimal hypersurface $\Sigma$ can be viewed in two ways.  One is that $\Sigma$ is $f$-minimal if and only if $\Sigma$ is a critical point of the weighted volume functional $e^{-f}d\sigma$,  where $d\sigma$ is the volume element of $\Sigma$. The other one is that $\Sigma$ is $f$-minimal if and only if $\Sigma$ is minimal in the new conformal metric $\tilde{g}=e^{-\frac{2f}{n}}\bar {g}$ (see Section \ref{notation} and Appendix).  $f$-minimal hypersurfaces have been studied before, even  more general stationary hypersurfaces for parametric elliptic functionals,   see for instance the work of White \cite{W} and Colding-Minicozzi \cite{CM5}.

 We prove the following compactness result:
\begin{theorem}\label{prin-intro}
Let $(M^{3},\bar{g},e^{-f}d\mu)$ be a complete smooth metric measure space and $\overline{\textrm{Ric}}_{f}\geq k$, where $k$ is a positive constant. Given an integer $g\geq0$ and a constant $V>0$, the space $S_{g, V}$ of smooth complete embedded $f$-minimal surfaces $\Sigma\subset M$ with 
\begin{itemize}
\item genus at most $g$,
\item $\partial\Sigma=\emptyset$, 
\item  $\int_{\Sigma}e^{-f}d\sigma\leq V$
\end{itemize}
is compact in the $C^m$ topology, for any $m\geq 2$.
Namely, any sequence of $S_{g, V}$ has a subsequence that converges in the $C^m$ topology on compact subsets to a surface in $S_{D, g}$, for any $m\geq 2$.
\end{theorem}

Since  the existence of  the uniform scale-invariant area bound is equivalent to the existence of  the uniform bound of the weighted area for self-shrinkers (see Remark \ref{general} in Section \ref{compactness2}),  Theorem \ref{prin-intro} implies Theorem \ref{cm-compact}. Also, in \cite{CMZ1}, we will apply Theorem \ref{prin-intro} to obtain a compactness theorem for the space of closed embedded $f$-minimal surfaces with the upper bounds of genus and diameter.

To prove Theorem \ref{prin-intro},  we  need to prove a non-existence  result on $L_f$-stable $f$-minimal hypersurfaces, which is of independent interest.
\begin{theorem}\label{stable-intro}
Let $(M^{n+1}, \bar{g},e^{-f}d\mu)$ be a complete smooth metric measure space with $\overline{\textrm{Ric}}_{f}\geq k$, where $k$ is positive constant. Then there is no   complete two-sided $L_{f}$-stable  $f$-minimal hypersurface  $\Sigma$ immersed in $(M, \bar{g})$ without
boundary and with finite weighted volume (i.e. $\int_\Sigma e^{-f}d\sigma<\infty$), where $d\sigma$ denotes the volume element on $\Sigma$ determined by the induced metric from $(M,\overline{g})$.
\end{theorem}

Here we explain briefly the meaning of $L_f$ stability.  For an $f$-minimal hypersurface $\Sigma$,  $L_f$ operator is $$L_f=\Delta_f+|A|^2+\overline{Ric}_f(\nu,\nu),$$ where $\Delta_f =\Delta-\langle\nabla f,\nabla\cdot\rangle$ is the weighted Laplacian on $\Sigma$. 

Especially, for self-shrinkers, it is so-called $L$ operator:
$$L=\Delta-|A|^2-\frac12\langle x,\nabla \cdot\rangle +\frac12.$$

$L_f$-stability of $\Sigma$ means that its weighted volume $\int_{\Sigma}e^{-f}d\sigma$ is locally minimal, that is,  the second variation of its weighted volume is nonnegative for all compactly supported normal variation. We leave more details about  the definition of  $L_f$-stability  and some of its properties in Section \ref{notation} and Appendix.

For self-shrinkers in $\mathbb{R}^{n+1}$, Colding-Minicozzi \cite{CM4} proved that
\begin{theorem}\cite{CM4}\label{cm} There are no $L$-stable smooth complete self-shrinkers without boundary and with polynomial volume growth in $\mathbb{R}^{n+1}$.
\end{theorem}
Since  the first and third authors  \cite{XZ} of the present paper proved  that for self-shrinkers, properness, the polynomial volume growth, and finite weighted volume are equivalent. Hence  Theorem \ref{stable-intro} implies Theorem \ref{cm}.  

In this paper,  we discuss  the relation among the properness, polynomial volume growth and  finite weighted volume of  $f$-minimal submanifolds (Propositions \ref{prop-imer1}, \ref{poli-grow} and \ref{finit-weig}). We obtain their  equivalence when the ambient space $(M,\bar{g},f)$ is a shrinking gradient Ricci solitons, i.e. $\overline{\textrm{Ric}}+\overline{\nabla}^2f=\frac12 \bar{g}$ with the condition $|\overline{\nabla}f|^2\leq f$ (Corollary \ref{equiv}).

The rest of this paper is organized as follows: In Section \ref{notation}  some definitions, notations and facts are given as a preliminary; In Section \ref{relation}  we prove Propositions \ref{prop-imer1}, \ref{poli-grow} and \ref{finit-weig};  In Section \ref{sec-stable} we prove Theorem \ref{stable-intro}; In Section \ref{compactness2} we prove Theorem \ref{prin-intro};  In Appendix we discuss some properties of $L_f$-stability for $f$-minimal submanifolds.

\section{Preliminaries}\label{notation}

In general, a smooth metric measure space, denoted by $(M^m, \bar{g}, e^{-f}d\mu)$, is an $m$-dimensional Riemannian manifold $(M^m, \bar{g})$ together with a weighted volume form $e^{-f}d\mu$ on $M$, where $f$ is a smooth function on $M$ and $d\mu$ the volume element induced by the metric $\bar{g}$.   In this paper, unless otherwise specified, we denote by a bar all quantities on $(M, \bar{g})$, for instance by $\overline{\nabla}$ and $\overline{\textrm{Ric}}$, 
 the Levi-Civita connection  and  the  Ricci curvature tensor of $(M, \bar{g})$ respectively.
For $(M, \bar{g}, e^{-f}d\mu)$, an  important and natural tensor is the $\infty$-Bakry-\'Emery Ricci curvature tensor $\overline{\text{Ric}}_{f}$ (for simplicity, Bakry-\'Emery Ricci curvature), which  is defined by
$$
\overline{\text{Ric}}_{f }:=\overline{\text{Ric}}+\overline{\nabla}^{2}f,
$$
where  $\overline{\nabla}^{2}f$ is the Hessian of $f$ on $M$.
If $f$ is constant,   $\overline{\textrm{Ric}}_{f}$ is the  Ricci curvature  $\overline{\textrm{Ric}}$ on $M$ respectively.

A Riemannian manifold with Bakry-\'Emery Ricci curvature bounded below by a positive constant has some properties similar to a Riemannian manifold with Ricci curvature bounded below by a positive constant. For instance, see the work of Wei-Wylie \cite{WW}, Munteanu-Wang \cites{MW, MW1} and the references therein. In this paper,   we will use  the following proposition by Morgan \cite{M} (see also its proof in  \cite{WW} ).
\begin{prop}\label{group} If a complete smooth metric measure space  $(M, \bar{g}, e^{-f}du)$ has $\overline{\text{Ric}}_{f}\geq k$, where $k$ is a positive constant, then $M$ has finite weighted volume (i.e. $\int_Me^{-f}d\mu<\infty$) and finite fundamental group.
\end{prop}

Now,  let $i: \Sigma^{n}\to M^m, n<m,$ be an $n$-dimensional smooth immersion.
Then $i: (\Sigma^{n}; i^*\bar{g}) \to (M^m, \bar{g})$ is an isometric immersion with the induced   metric $i^*\bar{g}$. For simplicity, we still  denote $i^*\bar{g}$  by $\bar{g}$ whenever  there is no confusion. We will denote for instance by ${\nabla}$,  $\textrm{Ric}$, $\Delta$ and $d\sigma$, 
 the Levi-Civita connection, the  Ricci curvature tensor,  the Laplacian, and the volume element of $(\Sigma, \bar{g})$ respectively.

The function $f$ induces a weighted measure  $e^{-f}d\sigma$ on $\Sigma$. Thus we have   an induced smooth metric measure space $(\Sigma^{n}, \bar{g}, e^{-f}d\sigma)$. 
 
 The associated weighted Laplacian ${\Delta}_{f}$ on  $(\Sigma, \bar{g})$ is defined  by
$$
{\Delta}_{f}u:={\Delta} u-\langle{\nabla} f,{\nabla} u\rangle.
$$
 The second order operator ${\Delta}_f$  is a self-adjoint operator on the space of square integrable functions on $\Sigma$ with respect to the measure $e^{-f}d\sigma$ (however the Laplacian operator in general has no this properties).

The second fundamental form $A$ of $(\Sigma,\bar{g})$ is defined by $$A(X,Y)=(\overline{\nabla}_XY)^{\perp}, \quad X, Y\in T_p \Sigma,p\in \Sigma,$$
where  $\perp$ denotes the projection to the normal bundle of $\Sigma$. 
The mean curvature vector ${\bf H}$ of $\Sigma$ is defined by ${\bf H}=\text{tr}A=\sum_{i=1}^{n}(\overline{\nabla}_{e_{i}}e_{i})^{\bot}$.

 \begin{definition}
The weighted mean curvature vector of $\Sigma$ with respect to  the metric $\bar{g}$ is defined by
\begin{equation}{\bf H}_{f}={\bf H}+(\bar{\nabla} f)^{\bot}.
\end{equation}
\end{definition}

The immersed submanifold $(\Sigma,\bar{g})$ is called $f$-minimal if its weighted mean curvature vector ${\bf H}_f$ vanishes identically, or equivalently if its  mean curvature vector  satisfies
\begin{equation}
{\bf H}=-(\bar{\nabla}f)^{\bot}.
\end{equation}

\begin{definition}The weighted volume  of $(\Sigma, \bar{g})$ is defined by
\begin{equation}V_f(\Sigma):=\int_\Sigma e^{-f}d\sigma.
\end{equation}

\end{definition}

It is well known  that  $\Sigma$ is $f$-minimal if and only if $\Sigma$ is a critical point of the weighted volume functional. Namely,  it  holds that
\begin{prop}\label{first-var}
If $T$ is a compactly supported variational field on $\Sigma$, then the first variation formula of the weighted volume of $(\Sigma,\bar{g})$ is given by
\begin{equation}\label{prop-fir-var}
\frac{d}{dt}V_{f}(\Sigma_{t})\biggr|_{t=0}=-\int_{\Sigma}\langle T^{\bot}, {\bf H}_{f}\rangle_{\bar{g}} e^{-f}d\sigma.
\end{equation}
\end{prop}

On the other hand,  an $f$-minimal submanifold can be viewed as a minimal submanifold under a conformal metric. Precisely,  define the new metric  $\tilde{g}=e^{-\frac{2}{n}f}\bar{g}$ on $M$, which is  conformal to $\bar{g}$.  Then the immersion $i: \Sigma\to M$ induces a metric $i^*\tilde{g}$ on $\Sigma$ from $(M, \tilde{g})$.   
In the following, $i^*\tilde{g}$ is still denoted by $\tilde{g}$ for simplicity. The volume of $(\Sigma, \tilde{g})$   is
\begin{equation}\label{two-vol}
\tilde V(\Sigma):=\int_\Sigma d\tilde{\sigma}=\int_\Sigma e^{-f}d\sigma=\text{V}_f(\Sigma).
\end{equation}
\noindent Hence Proposition \ref{first-var} and (\ref{two-vol}) imply that 
\begin{equation}\label{fir-variation}
\int_{\Sigma}\langle T^{\bot}, \tilde{{\bf H}}\rangle_{\tilde{g}} d\tilde{\sigma}=\int_{\Sigma}\langle T^{\bot}, {\bf H}_{f}\rangle_{\bar{g}} e^{-f}d\sigma,
\end{equation}
where $d\tilde{\sigma}=e^{-f}d\sigma$ and $ \tilde{{\bf H}}$ denote the volume element and the mean curvature vector of $\Sigma$ with respect to the conformal metric $\tilde{g}$ respectively.

Identity (\ref{fir-variation}) implies that $ \tilde{{\bf H}}=e^{\frac{2f}{n}}{\bf H}_f$ and  $(\Sigma, \bar{g})$ is $f$-minimal in $(M,\bar{g})$ if and only if $(\Sigma, \tilde{g})$ is minimal in $(M,\tilde{g})$.

Now suppose that    $\Sigma^n$ is a hypersurface immersed in $M^{n+1}$. Let $p\in \Sigma$ and  $\nu$ a unit normal at $p$. The second fundamental form $A$ and the mean curvature $H$ of $(\Sigma,\bar{g})$ are   as follows:  $$A: T_p\Sigma\to T_p\Sigma, A(X)=\overline{\nabla}_X\nu, X\in T_p\Sigma,$$  $$H=\text{tr}A=-\displaystyle\sum_{i=1}^n\langle \overline{\nabla}_{e_i}e_i,\nu\rangle.$$ Hence the mean curvature vector $\bf H$ of $(\Sigma,\bar{g})$ satisfies ${\bf H}=-H\nu.$  Define the weighted mean curvature $H_f$ of $(\Sigma,\bar{g})$ by ${\bf H}_f:=-H_f\nu$. Then $$H_f=H-\langle \overline{\nabla} f,\nu\rangle.$$

\begin{definition} A hypersurface $\Sigma$ immersed in $(M^{n+1}, \bar{g}, e^{-f}d\mu)$ with the induced metric $\bar{g}$  is called an $f$-minimal hypersurface if it satisfies
\begin{equation} H=\langle \overline{\nabla} f,\nu\rangle.
\end{equation}
\end{definition}

For  a hypersurface $(\Sigma,\bar{g})$, the $L_f$ operator is defined by $$L_f:=\Delta_f+|A|^2+\overline{\text{Ric}}_f(\nu,\nu),$$ where $|A|^2$ denotes the square of the norm of the second fundamental form $A$ of $\Sigma$. 

The $L_f$-stability of $\Sigma$ is defined as follows:
\begin{definition}A two-sided $f$-minimal hypersurface  $\Sigma$ is said to be $L_{f}$-stable if  for any  compactly supported  smooth function $\phi\in C_o^{\infty}(\Sigma)$, it holds that
\begin{equation}
-\int_{\Sigma}\phi L_f\phi e^{-f}d\sigma=\int_{\Sigma}[|\nabla\phi|^2-(|A|^{2}+\overline{\textrm{Ric}}_{f}(\nu,\nu))\phi^2]e^{-f}d\sigma\geq 0.
\end{equation}
\end{definition}
It is known that an   $f$-minimal hypersurface  $(\Sigma,\bar{g})$ is $L_f$-stable if and only if  $(\Sigma, \tilde{g})$ is stable as a minimal surface with respect to the conformal metric $\tilde{g}=e^{-f}\bar{g}$. See more details in Appendix of this paper.

In this paper,  for closed hypersurfaces, we choose $\nu$ to be  the outer unit normal.

\section{Properness, polynomial volume growth and  finit weighted volume of  $f$-minimal hypersurfaces}\label{relation}

In \cite{XZ}, the first and  third authors of the present paper  proved that the finite weighted volume of a self-shrinker $\Sigma^n$ immersed in $\mathbb{R}^{m}$ implies it is properly immersed.  In \cite{DX}, Ding-Xin proved that a properly immersed self-shrinker must have the Euclidean volume growth.  Combining these two results, it was proved \cite{XZ} that for immersed self-shrinkers, properness, polynomial volume growth and  finite weighted volume are equivalent.

In this section we  study the relation among the properness, polynomial volume growth and  finite weighted volume of  $f$-minimal submaifolds, some of them  will be used later in this paper.  

 If $\Sigma$ is an $n$-dimensional submanifold in a complete manifold $M^m,  n<m$,  $\Sigma$ is said to have polynomial volume growth  if, for a $p\in M$ fixed, there exist constants $C$ and $d$ so that for all $r\geq 1$,
\begin{equation} \text{Vol} ({B}^M_r( p)\cap \Sigma)\leq Cr^d \label{P-1-1},
\end{equation}
where ${B}^M_r( p)$ is the extrinsic ball of radius $r$ centered at $p$,  $ \text{Vol} ({B}^M_r( p)$ denotes the volume of ${B}^M_r( p)\cap \Sigma$.
When $d=n$ in (\ref{P-1-1}), $\Sigma$ is said to be of  Euclidean volume growth.

Before proving the following Proposition \ref{prop-imer1}, we recall  an estimate implied by  the Hessian comparison theorem  (cf, for instance,~\cite{CM4} {Lemma 7.1}).
\begin{lemma}\label{prop-imer}
Let $(M,\bar{g})$ be a complete Riemannian manifold with bounded geometry, that is, $M$ has  sectional curvature bounded by $k$ ($|K_M|\leq k$), and injectivity radius bounded below by $i_{0}>0$. Then the distance function $r(x)$ satisfies
\begin{equation*}
\bigr|\overline{\nabla}^{2}r(V,V)-\frac{1}{r}\bigr|V-\langle V,\overline{\nabla}r\rangle\overline{\nabla}r\bigr|^{2}\bigr|\leq\sqrt{k},
\end{equation*}
for $r<\min\{i_{0},\frac{1}{\sqrt{k}}\}$ and any unit vector $V\in T_x\Sigma$.
\end{lemma}
Using this estimate we will prove that
\begin{prop}\label{prop-imer1} Let $\Sigma^{n}$ be a complete noncompact $f$-minimal submanifold immersed in a complete Riemannian manifold  $M^m$. If $\Sigma$ has finite weighted area, then  $\Sigma$ is properly immersed in $M$.
\end{prop}
\noindent\textbf{Proof}. We argue by contradiction. Since  the argument is local, we may assume that $(M,g)$ has bounded geometry. Suppose that  $\Sigma$ is not properly immersed.  Then there exist a number $2R<\min\{i_{0},\frac{1}{\sqrt{k}}\}$ and $o\in M$ so that  $\overline{B}^{M}_{R}(o)\cap\Sigma$ is not compact in $\Sigma$, where $B^{M}_{R}(o)$ denotes an extrinsic ball of  radius $R$ centered at $o$. Then for any $a>0$, there is a sequence $\{p_{k}\}$ of points in $B^{M}_{R}(o)\cap\Sigma$ with $\textrm{dist}_{\Sigma}(p_{k},p_{j})\geq a>0$ for any $k\neq j$. So $B_{\frac{a}{2}}^{\Sigma}(p_{k})\cap B_{\frac{a}{2}}^{\Sigma}(p_{j})=\emptyset$ for any $k\neq j$. Choose $a<2R$. Then $ B_{\frac{a}{2}}^{\Sigma}(p_{j})\subset B^{M}_{2R}(o)$.  If $p\in B_{\frac{a}{2}}^{\Sigma}(p_{j})$, the extrinsic distance function $r_{j}(p)=\textrm{dist}_{M}(p,p_{j})$ from $p_{j}$ satisfies,
\begin{align*}
\Delta r_{j}&=\sum_{i=1}^{n}\overline{\nabla}^{2}r_{j}(e_{i},e_{i})+\langle{\bf H},\overline{\nabla}r_{j}\rangle\\
&\geq\frac{n}{r_{j}}-\frac{1}{r_{j}}|\nabla r_{j}|^{2}-n\sqrt{k}-\langle\overline{\nabla}f^{\bot},\overline{\nabla}r_{j}\rangle\\
&\geq\frac{n}{r_{j}}-\frac{1}{r_{j}}|\nabla r_{j}|^{2}-c,
\end{align*}
where $c=n\sqrt{k}+\sup_{B_{2R}^{M}(0)}|\overline{\nabla}f|$. Hence
\begin{equation*}
\Delta r^{2}_{j}\geq 2n-2cr_{j}.
\end{equation*}
Choosing $a\leq\min\{\frac{n}{2c},2R\}$, we have for $0<\mu\leq\frac{a}{2}$
 \begin{align} \label{P-4}
 \int_{B^{\Sigma}_{\mu}( p_j)} (2n-2cr_j)d\sigma&\le  \int_{B^{\Sigma}_{\mu}( p_j)}     \Delta_\Sigma r_j^2d\sigma\\
 &=  \int_{\partial B^{\Sigma}_{\mu}( p_j)}    \langle \nabla^{\Sigma} r_j^2, \nu\rangle d\sigma\nonumber\\
 &\le 2\mu A(\mu),\nonumber
 \end{align}
 where $\nu$ denotes the  outward unit normal vector of  $\partial B^{\Sigma}_{\mu}( p_j)$ and $A(\mu)$ denotes the area of $\partial B^{\Sigma}_{\mu}( p_j)$. Using co-area formula in (\ref{P-4}), we have
\begin{equation}\int_0^\mu (n- cs)A(s)ds\le \int_0^{\mu}\int_{d_{\Sigma}(p,p_j)=s}(n-cr_j)d\sigma\le\mu A(\mu). \label{P-5}
\end{equation}
This implies
\begin{equation*}
(n-c\mu)V(\mu)\leq V'(\mu),
\end{equation*}
where $V(\mu)$ denotes the volume of $B^{\Sigma}_{\mu}(p_{j})$.  So
\begin{equation}\frac{V'(\mu)}{V(\mu)}\geq \frac{n}{\mu}- c \label{P-6}
\end{equation}
Integrating (\ref{P-6}) from $\epsilon>0$ to $\mu$,  we have
\begin{equation}\frac{V(\mu)}{V(\epsilon)}\geq  (\frac{\mu}{\epsilon})^ne^{-c(\mu-\epsilon)}\nonumber
\end{equation}
Since
$\displaystyle\lim_{s\rightarrow 0^+}\frac{V(s)}{s^{n}}=\omega_{n}$,
\begin{equation}\label{P-7}V(\mu)\ge \omega_{n}\mu^ne^{- c\mu}.
\end{equation}
Thus we conclude
\begin{equation*}
\int_{\Sigma}e^{-f}d\sigma\geq\sum_{j=1}^{\infty}\int_{B^{\Sigma}_{\frac{a}{2}}(p_{j})}e^{-f}d\sigma\geq \inf_{B_{2R}^{M}(o)}(e^{-f})\sum_{j=1}^{\infty}\int_{B^{\Sigma}_{\frac{a}{2}}(p_{j})}d\sigma=\infty.
\end{equation*}
This contradicts with the assumption of the finite weighted volume of $\Sigma$.

\qed

\begin{prop}\label{poli-grow}
Let $(M^{m},\bar{g},e^{-f}d\mu)$ be a complete smooth metric measure space with $\overline{\Ric}_{f}=k,$ where $k$ is a positive constant. Assume that $f$ is a convex function. If  $\Sigma^{n}$ is a  complete noncompact properly immersed $f$-minimal submanifold, then $\Sigma$ has finite weighted volume and Euclidean (hence polynomial) volume growth.
\end{prop}
\noindent\textbf{Proof}. Since  $(M,\bar{g},f)$ is a gradient shrinking Ricci soliton,  it is well-known that, by a scaling  of the metric $\bar{g}$ and a translating of $f$, still denoted by $\bar{g}$ and $f$, we may normalize the metric so that $k=\frac{1}{2}$ and  the following identities hold:
\begin{align*}
\bar{R}+|\overline{\nabla}f|^{2}-f&=0,\\
\bar{R}+\bar{\Delta}f&=\frac{m}{2},\\
\bar{R}&\geq0.
\end{align*}
From these equations, we have that
\begin{equation*}
\bar{\Delta}f-|\overline{\nabla}f|^{2}+f=\frac{m}{2},\qquad \textrm{and}\qquad |\overline{\nabla}f|^{2}\leq f.
\end{equation*}
It was  proved by Cao and the third author \cite{CZ} that there is a positive constant $c$ so that
\begin{equation}\label{potencial-1}
\frac{1}{4}(r(x)-c)^{2}\leq f(x)\leq\frac{1}{4}(r(x)+c)^{2}.
\end{equation}
for any $x\in M$ with $r(x)=\textrm{dist}_{M}(p,x)\geq r_{0}$, where $p$ is an fixed point in $M$ and $c,r_{0}$ are positive constants that depend only of $m$ and $f(p)$.

By (\ref{potencial-1}), we know that $f$ is a proper function on $M$.
Since $\Sigma$ is properly immersed in $M$ and $f$ is proper in $M$,  $f|_{\Sigma}$ is also a  proper smooth function on $\Sigma$. Note that with  the scaling  metric and translating $f$,  $\Sigma$ is still $f$-minimal.
Hence
\begin{align*}
\Delta f-|\nabla f|^{2}+f&=(\bar{\Delta}f-\sum_{\alpha=n+1}^{m}f_{\alpha\alpha}-|\overline{\nabla}f^{\bot}|^{2})-|\overline{\nabla} f^{\top}|^{2}+f\\
&=\bar{\Delta}f-|\overline{\nabla}f|^{2}+f-\sum_{\alpha=n+1}^{m}f_{\alpha\alpha}\\
&\leq\frac{m}{2}.
\end{align*}
Also we have
\begin{equation*}
|\nabla f|^{2}=|\overline{\nabla}f^{\top}|^{2}\leq|\overline{\nabla}f|^{2}\leq f.
\end{equation*}
By Theorem 1.1 of~\cite{XZ}, $\Sigma$ has finite weighted volume and the Euclidean volume growth of the sub-level set of $f$  with the respect to the scaling metric and the translating $f$, and  hence with the respect to the original metric and $f$. Moreover, by the estimate (\ref{potencial-1}), we have that $\Sigma$ has the Euclidean volume growth.

\qed

We prove the following
\begin{prop}\label{finit-weig}
Let $(M^{m},\bar{g},e^{-f}d\mu)$ be a complete smooth metric measure space with $\overline{\Ric}_{f}\geq k$, where $k$ is a positive constant. Assume that   $|\overline{\nabla}f|^{2}\leq  2kf$. If  $\Sigma^{n}$ is a complete submanifold (not necessarily $f$-minimal) with polynomial area growth, then $\Sigma$ has finite weighted volume.
\end{prop}
\noindent\textbf{Proof}. By a scaling of the metric, we may assume that $k= \frac{1}{2}$. The proof follows from the  estimate of  $f$.  Munteanu-Wang~\cite{MW1} extended the estimate (\ref{potencial-1})  to $(M^{m},\bar{g},e^{-f}d\mu)$ with $\overline{\Ric}_{f}\geq \frac12$ and  $|\overline{\nabla}f|^{2}\leq  f$.
Combining the assumption that $\Sigma$ has polynomial volume growth with the estimative~\eqref{potencial-1}, we have
\begin{align*}
\int_{\Sigma}e^{-f}d\sigma&=\int_{\Sigma\cap B^{M}_{r_{0}}(p)}e^{-f}d\sigma+\sum_{i=0}^{\infty}\int_{\Sigma\cap(B^{M}_{r_{0}+i+1}(p)\backslash B^{M}_{r_{0}+i}(p))}e^{-f}d\sigma,\\
&\leq C_{1}\text{Vol}(\Sigma\cap B^{M}_{r_{0}}(p))+C\sum_{i=0}^{\infty}e^{-\frac{1}{4}(r_{0}+i-c)^{2}}\text{Vol}(\Sigma\cap B_{r_{0}+i+1}^{M}(p)),\\
&\leq C\bigr[r_{0}^{d}+\sum_{i=0}^{\infty}e^{-\frac{1}{4}(r_{0}+i-c)^{2}}(r_{0}+i+1)^{d}\bigr],\\
&< \infty.
\end{align*}
\qed

By Propositions \ref{prop-imer1}, \ref{poli-grow} and \ref{finit-weig}, we have the following
\begin{cor}\label{equiv}Let $(M^{m},\bar{g},f)$ be a complete shrinking gradient Ricci soliton with $\overline{\Ric}_{f}\geq \frac12$. Assume that   $|\overline{\nabla}f|^{2}\leq  f$.  If $\Sigma$ is a complete $f$-minimal submanifold immersed in $M$, then  for  $\Sigma$ the properness, polynomial volume growth, and finite weighted volume are equivelent.
\end{cor}

\section{Non-existence of $L_f $ stabe $f$-minimal hypersurfaces}\label{sec-stable}

In this section, we prove Theorem \ref{stable-intro}, which is a key to prove the compactness theorem in Section \ref{compactness2}.
\begin{theorem}\label{stable}(Theorem \ref{stable-intro})
Let $(M, \bar{g},e^{-f}d\mu)$ be a complete smooth metric measure space with $\overline{\textrm{Ric}}_{f}\geq k$, where $k$ is a positive constant. Then there is no two-sided $L_{f}$-stable complete  $f$-minimal hypersurface  $\Sigma$ immersed in $(M,g)$ without boundary and with finite weighted volume  (i.e. $\int_{\Sigma}e^{-f}d\sigma<\infty$).
\end{theorem}
\noindent\textbf{Proof}.  We argue by contradiction. 
Suppose that $\Sigma$ is an $L_{f}$-stable complete  $f$-minimal hypersurface  immersed in $(M,g)$ without boundary and with finite weighted volume. Recall that a two-sided  hypersurface $\Sigma$ is  $L_{f}$-stable if the following inequality holds that,  for any compactly supported  smooth function  $\phi\in\mathcal{C}^{\infty}_{o}(\Sigma)$,
\begin{equation}\label{f-stable-ine}
\int_{\Sigma}\left[|\nabla\phi|^{2}-\bigr(|A|^{2}+\overline{\textrm{Ric}}_{f}(\nu,\nu)\bigr)\right]e^{-f}d\sigma\geq 0.
\end{equation}

Observe that any closed hypersurface cannot be $L_{f}$-stable. This is because that: the assumption $\overline{\textrm{Ric}}_{f}\geq k>0$ implies that (\ref{f-stable-ine}) cannot hold for $\phi\equiv c$ on $\Sigma$. Hence, $\Sigma$ must be noncompact.

Let $\eta$ be  a nonnegative smooth function on $[0,\infty)$ satisfying
$$
\eta(s)=\left\{
\begin{array}{ccc}
1&\textrm{if}&s\in[0,1)\\
0&\textrm{if}&s\in[2,\infty)
\end{array}
\right.
$$
and $|\eta'|\leq2$.

Fix a point $p\in \Sigma$  and  let $r(x)=\textrm{dist}_{\Sigma}(p,x)$ denote the (intrinsic) distance function in $\Sigma$. Define a sequence of  functions $\phi_{j}(x)=\eta(\frac{r(x)}{j})$, $j\geq 1$.   Then  $|\nabla\ph_j|^{2}\leq1$ for $j\geq2$.  Substituting $\phi_j$ for $\phi$ in (\ref{f-stable-ine}):
\begin{align*}
\int_{\Sigma}\biggr(|\nabla\phi_j|^{2}
&-(|A|^2+\overline{\textrm{Ric}}_{f}(\nu,\nu))\phi_j^{2}\biggr)e^{-f}d\sigma\\
&\leq\int_{\Sigma}
                               \bigr(|\nabla\phi_j|^{2}-k\phi_j^{2}\bigr)e^{-f}d\sigma\\
                                                 &=\int_{B_{2j}^{\Sigma}( p)\backslash B_{j}^{\Sigma}( p)}|\nabla\phi_j|^{2}e^{-f}d\sigma-\int_{B_{2j}^{\Sigma}( p)
                                                    }k\phi_j^{2}e^{-f}d\sigma\\
                                                 &\leq\int_{B_{2j}^{\Sigma}( p)\backslash B_{j}^{\Sigma}( p)}e^{-f}d\sigma-k\int_{B_{2j}^{\Sigma}( p)}e^{-f}d\sigma\\
&\leq \int_{B_{2j}^{\Sigma}( p)\backslash B_{j}^{\Sigma}( p)}e^{-f}d\sigma-k\int_{B_{2}^{\Sigma}( p)}e^{-f}d\sigma,
\end{align*}
where $B_{j}^{\Sigma}( p)$ is the intrinsic geodesic ball in $M$ of radius $j$ centered at $p$. Since $\Sigma$ has finite weighted volume, we have, when $j\rightarrow\infty$
$$
\int_{B_{2j}^{\Sigma}( p)\backslash B_{j}^{\Sigma}( p)}e^{-f}d\sigma\rightarrow0.
$$
Choosing $j$ large enough, we have that $\phi_{j}$ satisfies
\begin{equation*}
\int_{\Sigma}\bigr(|\nabla\phi_j|^{2}-(|A|^{2}+\overline{\textrm{Ric}}_{f}(\nu,\nu))\phi_j^{2}\bigr)e^{-f}d\sigma<-\frac{k}{2}\int_{B_{2}^{\Sigma}( p)}e^{-f}d\sigma<0.
\end{equation*}
This contradicts that $\Sigma$ is $L_f$-stable.

\qed

\section{Compactness of complete $f$-minimal surfaces}\label{compactness2}

Before proving Theorem \ref{prin-intro}, we give some facts.

Wei-Wylie  (\cite{WW} Theorem 7.3) used mean curvature comparison theorem to give a distance estimate for two compact hypersurfaces  $\Sigma_1$ and $\Sigma_2$ in a smooth metric measure space $(M,\bar{g}, e^{-f}d\mu)$ with $\overline{\text{Ric}}_f \geq k$,  where $k$ is a positive constant.  Observe that  for two complete    hypersurfaces $\Sigma_1$ and $\Sigma_2$ if at least one of them is compact, there is a minimal geodesic joining  $\Sigma_1$ and $\Sigma_2$ and realizing their distance. Hence the  proof of Theorem 7.3 \cite{WW} can be applied to obtain the following

\begin{prop}\label{distance1}Let  $(M,\bar{g}, e^{-f}d\mu)$ be  an $(n+1)$-dimensional  smooth metric measure space with $\overline{\text{Ric}}_f \geq k$,  where $k$ is a positive constant.  If $\Sigma_1$ and $\Sigma_2$ are two complete  immersed    hypersurfaces, at least one of them is compact, then the distance $d(\Sigma_1,\Sigma_2)$ satisfies
\begin{equation}\label{distance}
d(\Sigma_1,\Sigma_2)\leq \frac{1}{k}(\max_{x\in \Sigma_1}|H^{\Sigma_1}_f( x)|+\max_{x\in\Sigma_2}|H_f^{\Sigma_2}(x)|).
\end{equation}
\end{prop}

\begin{cor} \label{inter}Let  $(M,\bar{g}, e^{-f}d\mu)$ be as in Proposition \ref{distance1}.  Then there is a  closed ball $B^M$ of $M$ satisfies that any complete immersed $ f$-minimal hypersurface $\Sigma$ must intersect it.
\end{cor}
\noindent \text{Proof}.  Fix $p\in M$ and a geodesic sphere $B^M_r( p)$ of $M$.  By Proposition \ref{distance1}, 
$$d(B^M_r( p), \Sigma)\leq \frac{1}{2k}\max_{x}|H^{B^M_r( p)}_f( x)|= C,$$
where $C$ is independent  on $\Sigma$.
Therefore there is a closed ball $B^M$ of $M$ with radius big enough so that any $\Sigma$ must intersect it.

\qed

We need  the following fact:

\begin{prop} \label{orient} Let $M$ be a simply connected Riemannian manifold. If an $f$-minimal hypersurface $\Sigma$ is complete, not necessarily connected, properly embedded, and has no boundary, then every component of $\Sigma$ separates $M$ into two components and  thus is two-sided. Therefore $\Sigma$ has a  globally defined unit normal.
\end{prop}

\noindent \textbf{Proof}.  Suppose $\Sigma_j$ is a component of $\Sigma$. By contrary, if $M\backslash \Sigma_j$ has one component. Since $\Sigma$ is a properly embedded $f$-minimal hypersurface, for any $p\in \Sigma_j$,  there is a neighborhood $W$ of $p$ in $M$ so that $W\cap \Sigma_j=W\cap\Sigma$ only has one piece (i.e. it is a graph above a connected domain in the tangent plane of $p$).  Thus we have a simply closed curve $\gamma$ passing $p$, transversal to $\Sigma_j$ at $p$, and $\Sigma_j\cap\gamma=p$.  Since $M$ is simply connected, we have a disk $D$ with the boundary $\gamma$.  Again since $\Sigma$ is proper,  the intersection of  $\Sigma_j$ with $\partial D=\gamma$ cannot be one point, which is a contradiction.

\qed

Combining  Proposition \ref{prop-imer1} in Section \ref{relation} with Proposition \ref{orient}, we obtain 

\begin{prop} \label{orient-1}Let $(M, \bar{g}, e^{-f}d\mu)$ be a simply connected complete smooth measure space.  If a complete $f$-minimal hypersurface has finite weighted volume, then every component of $\Sigma$ separates $M$ into two components and  thus is two-sided. Therefore $\Sigma$ has a  globally defined unit normal.

\end{prop}

We will take the same approach as in Colding-Minicozzi's paper  \cite{CM1} to prove Theorem  \ref{prin-intro}, a smooth compactness theorem for complete $f$-minimal surfaces. First we recall  a well known local singular compactness theorem for embedded minimal surfaces in a Riemannian $3$-manifold.
\begin{prop}\label{sing}  (cf \cite{CM1} Proposition 2.1)
Given a point $p$ in a Riemannian 3-manifold $M$,  there exists a $R>0$ such that : Let  $\Sigma_{j}$ be  embedded minimal surfaces in $B_{2R}(p)\subset M$ with $\partial\Sigma_{j}\subset\partial B_{2R}(p )$. If each $\Sigma_{j}$ has area at most $V$ and genus at most $g$ for some fixed $V,g$, then there is a finite collection of points $x_{k}$, a smooth embedded minimal surface $\Sigma\subset B_{R}(p )$ with $\partial\Sigma\subset\partial B_{R}(p )$ and a subsequence of  $\{\Sigma_{j}\}$ that converges in $B_{R}(p )$ (with finite multiplicity) to $\Sigma$ away from the set $\{x_{k}\}$.
\end{prop}

It is known that $\Sigma$ is $f$-minimal with respect to metric $\bar{g}$ if and only if $\Sigma$  is minimal with the conformal metric $\tilde{g}=e^{-f}g$ (see Appendix).  Using this fact and 
applying Proposition \ref{sing},  we may prove a global singular compactness theorem for $f$-minimal surfaces.
\begin{prop}\label{sing1} Let  $M$ be a complete $3$-manifold and $(M,\bar{g}, e^{-f}d\mu)$ a smooth metric measure  space.
Suppose that $\Sigma_{i}\subset M$ is a sequence of smooth complete embedded $f$-minimal surfaces with genus at most $g$,  without boundary, and with weighted area at most $V$, i.e.
\begin{equation}\label{comp1}
\int_{\Sigma_{i}}e^{-f}d\sigma\leq V<\infty.
\end{equation}
Then there is a subsequence,  still denoted by $\Sigma_{i}$, a smooth embedded complete non-trivial $f$-minimal surface $\Sigma\subset M$ without boundary, and a locally finite collection of points $\mathcal{S}\subset \Sigma$ so that $\Sigma_{i}$ converges smoothly (possibly with multiplicity) to $\Sigma$ off of $\mathcal{S}$. Moreover, $\Sigma$ satisfies $\int_{\Sigma}e^{-f}d\sigma\leq V$ and is properly embedded.
\end{prop}
\noindent Here a set $\mathcal{S}\subset M$ is said to be locally finite if $B_{R}(p)\cap\mathcal{S}$ is finite for every $p\in M$ and for all $R>0$.
\bigskip

\noindent\textbf{Proof}.  Consider the conformal metric $\tilde{g}=e^{-f}\bar{g}$ on $M$.  For a point  $p\in M$,  let $\tilde{B}_{2R}( p)\subset M$ denote  the ball in $(M, \tilde{g})$ of radius $2R$ centered at $p$.  Then the area of $\tilde{B}_{2R}(p )\cap\Sigma$   satisfies
\begin{equation}
\widetilde{\text{Area}}(\tilde{B}_{2R}(p )\cap\Sigma)\leq \int_{\Sigma}d\tilde{\sigma}=\int_{\Sigma}e^{-f}d\sigma\leq V.
\end{equation}
 Also, it is clear that the genus of $\tilde{B}_{2R}(p )\cap\Sigma_{j}$ remains at most $g$. Then by Proposition \ref{sing},  there exists a finite collection of points $x_{k}$, a smooth embedded minimal surface $\Sigma\subset \tilde{B}_{R}(p )$, with $\partial\Sigma\subset\partial \tilde{B}_{R}$ and a subsequence of  $\{\Sigma_{j}\}$ that converges in $\tilde{B}_{R}(p )$ (with finite multiplicity) to $\Sigma$ away from the set $\{x_{k}\}$.

Let $\{\tilde{B}_{R_i}(p_i)\}$ be  a countable   cover of  $(M,\tilde{g})$.  On each $\tilde{B}_{2R_i}(p_i)$, applying the previous local convergence, and then passing to a diagonal subsequence, we obtain that there is a subsequence of $\Sigma_i$,   still  denoted by  $\Sigma_i$,  a smooth embedded minimal surface $\Sigma$ (with respect to the metric $\tilde{g}$)  without boundary,  and a locally finite collection of points $\mathcal{S}\subset \Sigma$ so that $\Sigma_{i}$ converges smoothly (possibly with multiplicity) to $\Sigma$ off of $\mathcal{S}$. 

Since $\Sigma$ has no boundary,  it is complete in the original metric $\bar{g}$. Thus we obtain that the smooth convergence of the subsequence to the smooth embedded complete $f$-minimal surface $\Sigma$ off of $\mathcal{S}$.

By Corollary \ref{inter}, $\Sigma$ is nontrivial.

The convergence of $\Sigma_i$ to $\Sigma$ and (\ref{comp1}) imply $\int_{\Sigma}e^{-f}d\sigma\leq V$. By Proposition \ref{prop-imer1}, $\Sigma$ is properly embedded.

\qed

We need  to show that the convergence is smooth across the points in $\mathcal{S}$. To prove it, we need the following 
\begin{prop}\label{sing2} Assume that  the ambient manifold $M$ in Proposition \ref{sing1} is simply connected.
If the convergence of the sequence $\{\Sigma_{i}\}$  has the multiplicity greater than one, then $\Sigma$ is $L_{f}$-stable.
\end{prop}
\noindent\textbf{Proof}.  By Proposition \ref{orient-1}, we know that $\Sigma_i$ and $\Sigma$ are orientable. We may have two ways to prove the proposition. The first is to use the known fact on minimal surfaces. It is known  that  (cf \cite{CM4} Appendix A)  if the multiplicity of the  convergence of a sequence of embedded orientable minimal surfaces in a simply connected $3$-manifold is not one,  then the limit minimal surface is stable.  Under the conformal metric $\tilde{g}$,  a sequence $\{\Sigma_i\}$ of minimal surfaces converges to a smooth embedded orientable minimal surface $\Sigma$ and thus $\Sigma$ is stable.  Besides,  the conclusion that $\Sigma$ is stable with respect to the conformal metric $\tilde{g}$ is equivalent to that   $\Sigma$ is $L_f$-stable under the original metric $\bar{g}$ (see Appendix).

The second way is to prove directly. We may prove that  $L_{f}$ is the linearization of $f$-minimal equation by a similar proof to the one in \cites{CM3, CM2} Appendix A.  By arguing as in Proposition 3.2 in \cites{CM1, CM2},,  we can find a smooth positive function $u$ on $\Sigma$ satisfying
\begin{equation}
L_{f}u=0.
\end{equation}
This implies that $\Sigma$ is $L_{f}$-stable.

\qed

\noindent\textbf{Proof of theorem~\ref{prin-intro}}. By the assumption on $\overline{\textrm{Ric}}_f$ and Proposition \ref{group}, $M$ has finite fundamental group.  After passing to the universal  covering, we may assume that $M$ is simply connected.
Given a sequence of smooth complete embedded $f$-minimal surfaces $\{\Sigma_{i}\}$ with genus $g$, $\partial\Sigma_{i}=\emptyset$, and the weighted area at most $V$,  by Proposition \ref{sing1}, there is a  subsequence, still denoted by  $\{\Sigma_{i}\}$ so that  it converges in the topology of smooth convergence on compact subsets  to a smooth  embedded complete  $f$-minimal surface $\Sigma$ away from a locally finite set $\mathcal{S}\subset\Sigma$ (possibly with multiplicity). Moreover, the limit surface $\Sigma\subset M$ is complete, properly embedded, $\int_{\Sigma}e^{-f}d\sigma\leq V$, has no boundary and   has a well-defined unit normal $\nu$. We also have the equivalent convergence under the conformal metric $\bar{g}$.

If $\mathcal{S}$ is not empty,  Allard's regularity theorem implies that  the convergence has multiplicity greater than one. Then by Proposition \ref{sing2}, we conclude that $\Sigma$ is $L_{f}$-stable. But  Proposition~\ref{stable} says that  there is no such $\Sigma$.  This contradiction implies that $\mathcal{S}$ must be empty. We complete the proof of the theorem.

\qed

\begin{remark} \label{general} For self-shrinkers,  the condition that the scale-invariant uniform area bound exists (i.e. 
 there is a uniform bound $V_1$:  Area$(B_R(x_0)\cap\Sigma)\le V_1R^2$  for all $x_0\in \rr^3$ and $R>0$) implies that 
the uniform bound $V$ of weighted area  (i.e. $\int_{\Sigma}e^{-f}d\sigma<V$) exists (cf.  the proof of Proposition \ref{finit-weig}). The converse is also true by the conclusion that the entropy of a self-shrinker can be achieved by $F_{0,1}$  for self-shrinkers with polynomial volume growth (see Section 7 of \cite{CM3}). Therefore
Theorem \ref{prin-intro} generalizes the result of Colding-Minicozzi (Theorem \ref{cm-compact}) for self-shrinkers.

\end{remark}

\begin{remark} \label{compact1} Combining Theorem \ref{prin-intro} with the upper bound estimate of weighted area for closed  embedded $f$-minimal surfaces  of   fixed genus  in a complete $3$-manifold with $\overline{\textrm{Ric}}_{f}\geq k>0$, we may obtain the smooth compactness theorem for the space of  closed embedded $f$-minimal surfaces of fixed topological type and with  diameter bound. We discuss it in \cite{CMZ1}.

\end{remark}

\section{Appendix A}\label{appendix1}

In this Appendix,  we  discuss the $L_f$-stability properties of  $f$-submanifolds.  With the same notations as in Section \ref{notation}, let $(M^m, \bar{g})$ be an $m$-dimensional Riemannian manifold and  $i: \Sigma^{n}\to M^m, n<m,$ be  an immersion. Let  $\tilde{g}=e^{-\frac{2}{n}f}\bar{g}$ denote the new conformal metric on $M$.  Therefore  $i$  may induce two  isometric immersions of $\Sigma$:  $(\Sigma, \bar{g})\to (M,\bar{g})$ and $(\Sigma, \tilde{g})\to (M,\tilde{g})$ respectively.

When $(\Sigma, \tilde{g})$ is minimal, it is well known that the second variation of the  volume of $(\Sigma, \tilde{g})$ is given by
\begin{prop} (cf \cite{CM4}) Let $ (\Sigma, \tilde{g})$ be minimal submanifold in $(M, \tilde{g})$.
If $T$ is a normal compactly supported variational vector field on $\Sigma$ (that is, $T=T^\perp$), then the second variational formula of the volume $\tilde{V}$ of $(\Sigma, \tilde{g})$ is given by
\begin{equation}\label{sec-var}
\frac{d^{2}}{dt^{2}}\tilde{V}(\Sigma_{t})\biggr|_{t=0}=-\int_{\Sigma}\langle T,  J T\rangle_{\tilde{g}} d\tilde{\sigma},
\end{equation}
where  the stability operator (or Jacobi operator) $J$ is defined on a normal vector field $T$ to $\Sigma$ by
\begin{equation}JT=\Delta^{\perp}_{(\Sigma,\tilde{g})}T+tr_{(\Sigma,\tilde{g})} [\widetilde{Rm}(\cdot, T)\cdot ]^{\perp}+\tilde{B}(T).
\end{equation}
\noindent Here  $ \Delta^{\perp}_{(\Sigma,\tilde{g})}T=\displaystyle\sum_{i=1}^n(\nabla^{\perp}_{\tilde{e}_i}\nabla^{\perp}_{\tilde{e}_i}T-\nabla^{\perp}_{\nabla_{\tilde{e}_i}\tilde{e}_i}T)$ is the Laplacian determined by the normal connection $\nabla^{\perp}$ of $(\Sigma, \tilde{g})$, $\widetilde{Rm}$ is the curvature tensor on $(M,\tilde{g})$,  $\text{tr}_{(\Sigma,\tilde{g})} [\widetilde{Rm}(\cdot, T)\cdot ]^{\perp}=\displaystyle\sum_{i=1}^n[\widetilde{Rm}(\tilde{e}_i, T)\tilde{e}_i]^{\perp}$,
$\tilde{A}$ denotes the second fundamental form of $(\Sigma, \tilde{g})$,  $\tilde{B}(T)=\displaystyle\sum_{i,j=1}^n\langle \tilde{A}(\tilde{e}_i,\tilde{e}_j),T\rangle \tilde{A}(\tilde{e}_i,\tilde{e}_j)$,  and $\{\tilde{e}_i\}$, $i=1,\cdots, n$ is a  local orthonormal base of $(\Sigma,\tilde{g})$.
\end{prop}

Recall that the  weighted volume of $(\Sigma, \bar{g})$ is defined by
\begin{equation}V_f(\Sigma)=\int_\Sigma e^{-f}d\sigma.
\end{equation}

By a direct computation  similar to that of (\ref{sec-var}),  we may prove the second variational  formula of the weighted volume of $f$-minimal submanifold $(\Sigma, \bar{g})$.
\begin{definition}  For any normal vector field $T$ on $(\Sigma,\bar{g})$, the second order operator $ \Delta^{\perp}_f$ is defined by
\begin{align*} 
\Delta^{\perp}_fT&:=\Delta^{\perp}T-\text{tr}[\nabla f\otimes\nabla^{\perp}T(\cdot,\cdot)]\\
&=\displaystyle\sum_{i=1}^n(\nabla^{\perp}_{e_i}\nabla^{\perp}_{e_i}T-\nabla^{\perp}_{\nabla_{e_i}e_i}T)-\displaystyle\sum_{i=1}^n({e_i}f)(\nabla_{e_i}^{\perp}T)
\end{align*}

The  operator $L_f$ on $(\Sigma,\bar{g})$  is defined by,
\begin{equation}
L_fT=\Delta_f^{\perp}T+R(T)+B(T)+F(T).
\end{equation}
In the above,  $\nabla^{\perp}$ denotes the normal connection  of $(\Sigma, \tilde{g})$; $\{e_i\}$, $i=1,\ldots,n$ is a local orthonormal base of  $(\Sigma,\bar{g})$; $B(T)=\displaystyle\sum_{i,j=1}^n\langle{A}(e_i,e_j),T\rangle {A}(e_i,e_j)$,  where $A$ denotes the second fundamental form of $(\Sigma, \bar{g})$;

$R(T)=\text{tr}_{(\Sigma,\bar{g})} [\overline{Rm}(\cdot, T)\cdot ]^{\perp}=\displaystyle\sum_{i=1}^n[\overline{Rm}(e_i, T)e_i]^{\perp}$,
where $\overline{Rm}$ denotes the Riemannian curvature tensor of $(M,\bar{g})$;  and

 $F(T)=[\overline{\nabla}^2f(T)]^{\perp}=\displaystyle\sum_{\alpha=n+1}^m\overline{\nabla}^2f(T,e_{\alpha})e_{\alpha}$,
 where $\{e_{\alpha}\}$, $\alpha=n+1,\cdots,m$ is a local orthonormal normal vector field on $(\Sigma,\bar{g})$.
  \end{definition}
  \begin{prop} Let $(\Sigma, \bar{g})$ be an $f$-minimal in $(M, \bar{g})$.
If $T$ is a normal compactly supported variational vector field on $\Sigma$ (that is, $T=T^\perp$),  then the second variational of the weighted volume of $(\Sigma, \bar{g})$ is given by
\begin{equation}\label{Lfstable}
\frac{d^{2}}{dt^{2}}V_{f}(\Sigma_{t})\biggr|_{t=0}=-\int_{\Sigma}
\langle T,L_{f}T\rangle_{\bar{g}} e^{-f}d\sigma.
\end{equation}

\end{prop}

\noindent\textbf{Proof}.  
Let  $\psi(\cdot, t),  t\in(-\epsilon,\epsilon)$ be a compactly supported variation of $\Sigma$ so that $T=d\psi(\frac{\partial}{\partial t})$  is  the variational vector field,  $\Sigma_t=\psi(\Sigma, t), \Sigma_{0}=\Sigma$. Choose a normal coordinate system $\{x_{1},\ldots,x_{n}\}$  at a point $p\in \Sigma$. We can consider $\{x_{1},\ldots,x_{n},t\}$ to be a coordinate system of $\Sigma\times(-\epsilon,\epsilon)$ near the point $(p,0)$.   Denote $e_{i}=d\psi(\frac{\partial}{\partial x_{i}})$ for $i=1,\ldots,n$. The induced metric on $\Sigma_{t}$ from $(M,\bar{g})$ is  given for $g_{ij}=\langle e_{i},e_{j}\rangle$.  Hence $g_{ij}(p,0)=\delta_{ij}$ and $\nabla_{e_{i}}e_{j}(p,0)=0$. Denote by $d\sigma_{t}$ the volume element of $\Sigma_{t}$. Then $d\sigma_{t}=J(x,t)d\sigma_{0}$, where $d\sigma_{0}=d\sigma$  and the function $J(x,t)$ is given by
$$J(x,t)=\displaystyle\frac{\sqrt{G(x,t)}}{\sqrt{G(x,0)}},$$
with $G(x,t)=\textrm{det}(g_{ij}(x,t))$. Denote by $d(\sigma_{f})_{t}$ the weighted volume element of $\Sigma_{t}$. Then 
$d(\sigma_{f})_{t}=J_{f}(x,t)d\sigma_{0},$
where
$ J_{f}(x,t)=J(x,t)e^{-f(x,t)}, f(x,t)=f(\psi(x,t)).$

Since
$\frac{\partial J}{\partial t}=\sum_{i,j=1}^{n}g^{ij}\langle\overline{\nabla}_{e_{i}}T,e_{j}\rangle J,$
$\frac{\partial J_{f}}{\partial t}
=\bigr(\sum_{i,j=1}^{n}g^{ij}\langle\overline{\nabla}_{e_{i}}T,e_{j}\rangle-\langle\overline{\nabla} f,T\rangle \bigr)J_{f}.
$

Note that  $T$ is a  normal vector field. A direct computation gives, at $(p,0)$
 \begin{align*}
\frac{\partial^{2}J_{f}}{\partial^2 t}\biggr|_{t=0}&=\biggr[ -2\sum_{i,j=1}^{n}\langle A_{ij},T\rangle^2+\langle\bar{R}(e_{i},T)T,e_{i}\rangle\\
&\quad +\sum_{i=1}^{n}\langle\overline{\nabla}_{e_{i}}\overline{\nabla}_{T}T,e_{i}\rangle+\sum_{i=1}^{n}\langle
\overline{\nabla}_{e_{i}}T,\overline{\nabla}_{e_{i}}T
\rangle\\
&\quad  -\overline{\nabla}^{2}f(T,T)-\langle\overline{\nabla} f,\overline{\nabla}_{T}T\rangle\\
&\quad +\bigr(\sum_{i=1}^{n}\langle\overline{\nabla}_{e_{i}}T,e_{i}\rangle-\langle\overline{\nabla} f,T\rangle\bigr)\bigr(\sum_{j=1}^{n}\langle\overline{\nabla}_{e_{j}}T,e_{j}\rangle-\langle\overline{\nabla} f,T\rangle\bigr)\biggr]J_{f}.\\
\end{align*}
By
 \begin{align*}
 \sum_{i=1}^{n}\langle\overline{\nabla}_{e_{i}}T,\overline{\nabla}_{e_{i}}T\rangle&=\sum_{i,j=1}^{n}\langle\overline{\nabla}_{e_{i}}T,e_{j}
 \rangle^{2}+\sum_{i=1}^{n}
 \sum_{\alpha=n+1}^{m}\langle\overline{\nabla}_{e_{i}}T,e_{\alpha}\rangle^{2}\\
 &=\sum_{i,j=1}^{n}\langle A_{ij},T\rangle^{2}+\sum_{i=1}^{n}\langle\nabla^{\bot}_{e_{i}}T,\nabla^{\bot}_{e_{i}}T\rangle\\
 &=|\langle A(\cdot,\cdot),T\rangle|^{2}+|\nabla^{\bot}T|^{2},
 \end{align*}
and
  $ \sum_{i=1}^{n}\langle\overline{\nabla}_{e_{i}}\overline{\nabla}_{T}T,e_{i}\rangle=\textrm{div}(\overline{\nabla}_{T}T)^{\top}-\langle (\overline{\nabla}_{T}T)^{\bot}, {\bf H}\rangle, $
we have that, at $p$
\begin{align*}
\frac{\partial^{2}J_{f}}{\partial t^{2}}\biggr|_{t=0}&=\biggr[-|\langle A(\cdot,\cdot),T\rangle|^{2}-\sum_{i=1}^{n} \langle\bar{R}(e_{i},T)e_{i},T\rangle +|\nabla^{\bot}T|^{2}+\textrm{div}(\overline{\nabla}_{T}T)^{\top}\\
&\quad -\langle(\overline{\nabla}_{T}T)^{\bot},\vec{H}\rangle-\overline{\nabla}^{2}f(T,T)-\langle\overline{\nabla} f,\overline{\nabla}_{T}T\rangle+\langle T,{\bf H}_{f}\rangle^{2}\biggr]e^{-f}.
\end{align*}
 Using
$\textrm{div}\bigr(e^{-f}(\overline{\nabla}_{T}T)^{\top}\bigr)=e^{-f}\textrm{div}(\overline{\nabla}_{T}T)^{\top}-e^{-f}\langle(\overline{\nabla}_{T}T)^{\top},\nabla f\rangle$,
we have at $p$
\begin{align}\label{p}
&\frac{\partial^{2}J_{f}}{\partial t^{2}}\biggr[|_{t=0}\\
&=\bigg[|\nabla^{\bot}T|^{2}-|\langle A(\cdot,\cdot),T\rangle|^{2}-\sum_{i=1}^{n}\langle\bar{R}(e_{i},T)e_{i},T\rangle-\overline{\nabla}^{2}f(T,T)\nonumber\\
&\quad-\langle(\overline{\nabla}_{T}T)^{\bot},{\bf H}_{f}\rangle+\langle T,{\bf H}_{f}\rangle^{2}\biggr]e^{-f}+\textrm{div}\bigr(e^{-f}(\overline{\nabla}_{T}T)^{\top}\bigr).\nonumber
\end{align}
Observe that the right-hand side of (\ref{p}) is independent of the choice of coordinates. Hence (\ref{p}) holds on $\Sigma$.
By integrating (\ref{p}) and using the fact that $\Sigma$ is $f$-minimal (i.e., ${\bf H}_f=0$), we obtain
\begin{align*}
\frac{d^{2}}{dt^{2}}V_{f}(\Sigma_{t})\biggr|_{t=0}&=\int_{\Sigma}\bigr(|\nabla^{\bot}T|^{2}-|\langle A(\cdot,\cdot),T\rangle|^{2}-\langle{R}(T),T\rangle-\overline{\nabla}^{2}f(T,T)\bigr)e^{-f}d\sigma\\
&=-\int_{\Sigma}\langle T,\Delta^{\bot}_{f}T+{A}(T)+{R}(T)+F(T)\rangle e^{-f}d\sigma\\
&=-\int_{\Sigma}\langle T,L_{f}T\rangle e^{-f}d\sigma.
\end{align*}
 Substituting  $e^{-f}T$ for $T$ in the   identity  $\int_{\Sigma}|\nabla^{\bot}T|^{2}d\sigma=-\int_{\Sigma}\langle T,\Delta^{\bot}T\rangle d\sigma$, we have 
$$\int_{\Sigma}|\nabla^{\bot}T|^{2}e^{-f}d\sigma=-\int_{\Sigma}\langle T,\Delta^{\bot}_{f}T\rangle e^{-f}d\sigma.$$
Thus we have  the second variational formula of the weighted volume of $\Sigma$
\begin{align*}
\frac{d^{2}}{dt^{2}}V_{f}(\Sigma_{t})\biggr|_{t=0}
&=-\int_{\Sigma}\langle T,\Delta^{\bot}_{f}T+{A}(T)+{R}(T)+F(T)\rangle e^{-f}d\sigma\\
&=-\int_{\Sigma}\langle T,L_{f}T\rangle e^{-f}d\sigma.
\end{align*}
\qed

\begin{definition} An $f$-minimal submanifold $(\Sigma, \bar{g})$ is called $L_f$-stable if the second variational  of the weighted volume of $\Sigma$ given by (\ref{Lfstable}) is nonnegative for any  normal compactly supported variational vector field $T$ on $\Sigma$.
\end{definition}

Observe that for an $f$-minimal submanifold $\Sigma$ and its normal compactly supported variation,   it holds that $V_f(\Sigma_t)=\tilde{V}(\Sigma_t)$.  Then
\begin{equation}\label{eq1}
\frac{d^{2}}{dt^{2}}\tilde{V}(\Sigma_{t})\biggr|_{t=0}=\frac{d^{2}}{dt^{2}}V_{f}(\Sigma_{t})\biggr|_{t=0}.
\end{equation}
By (\ref{sec-var}),  (\ref{Lfstable}), and (\ref{eq1}),  we have
\begin{equation}\int_{\Sigma}\langle T,  J T\rangle_{\tilde{g}} e^{-f}d\tilde{\sigma}=\int_{\Sigma}
\langle T,L_{f}T\rangle_{\bar{g}} e^{-f}d\sigma.
\end{equation}
 This implies that
\begin{equation}\label{JL}
\int_{\Sigma}e^{-\frac{2f}{n}}\langle T,  J T\rangle_{\bar{g}} e^{-f}d\sigma=\int_{\Sigma}
\langle T,L_{f}T\rangle_{\bar{g}} e^{-f}d\sigma.
\end{equation}
By (\ref{JL}), the following equality holds.
\begin{cor} \label{j-lf}For any normal variation field $T$ on $\Sigma$,  $$JT=e^{\frac{2f}{n}}L_fT.$$
\end{cor}

The operator $L_f$ corresponds to  a  symmetric bilinear form $B_f(T,T)$  for the space of normal compactly supported field on $\Sigma$:
\begin{equation}
B_f(T,T):=-\int_{\Sigma}\langle T,  L_f T\rangle_{\bar{g}} e^{-f}d\sigma.
\end{equation}
 We define  the $L_{f}$-index, denoted by $L_{f}$-ind, of $(\Sigma, \bar{g})$  by the maximum of the dimensions of negative definite  subspaces of $B_f$.
Hence $(\Sigma, \bar{g})$ is $L_{f}$-stable if and only if its $L_{f}$-ind$=0$.

On the other hand, for minimal $(\Sigma, \tilde{g})$,  it is well known that  the stability operator $J$ also defines  a symmetric bilinear form $\tilde{B}(T,T)$,
\begin{equation}
\tilde{B}(T,T):=-\int_{\Sigma}\langle T,  J T\rangle_{\tilde{g}} d\tilde{\sigma}.
\end{equation}
There are also the concepts of  index and stability of $(\Sigma, \tilde{g})$.  In particular, $(\Sigma,\tilde{g})$ is stable if and only if the index $\text{ind}(\Sigma,\tilde{g})=0$. Since $B_f(T,T)=\tilde{B}(T,T)$, it holds that

\begin{prop}  $L_f$-ind of $(\Sigma, \bar{g})$ is equal to the index of $(\Sigma, \tilde{g})$. In particular,
 $(\Sigma, \bar{g})$ is $L_f$-stable if and only if $(\Sigma, \tilde{g})$ is stable in   $(M, \tilde{g})$.
\end{prop}

Now if $\Sigma$ is a two-sided hypersurface, that is, there  is a globally-defined unit normal $\nu$ on $(\Sigma,\bar{g})$. Take $T=\phi\nu$.  Then the second variation (\ref{Lfstable}) implies that
\begin{prop} Let $\Sigma$ be a two-sided $f$-minimal hypersurface in $(M^{n+1}, \bar{g})$.
If $\phi$ is a compactly supported smooth function on $\Sigma$, then the second variation of the weighted volume of $(\Sigma,  \bar{g})$ is given by
\begin{equation}
\frac{d^{2}}{dt^{2}}V_{f}(\Sigma_{t})\biggr|_{t=0}=-\int_{\Sigma}\phi L_{f}(\phi)e^{-f}d\sigma
\end{equation}
where  $\nu$ denotes the unit normal of $(\Sigma, \bar{g})$ and  the operator  $L_f$ is defined by $L_{f}=\Delta_{f}+|A|_{\bar{g}}^{2}+\overline{\textrm{Ric}}_{f}(\nu,\nu)$.
\end{prop}
\begin{definition}
The operator $L_{f}=\Delta_{f}+|A|_ {\bar{g}}^{2}+\overline{\textrm{Ric}}_{f}(\nu,\nu)$ is called $L_{f}$-stability operator of hypersurface$(\Sigma,  \bar{g})$.
\end{definition}
A bilinear form on space $C_o^{\infty}(\Sigma)$  of compactly supported smooth functions on $\Sigma$  is defined by
\begin{equation}\begin{split}
B_f(\phi, \phi):&=-\int_\Sigma \phi L_f\phi e^{-f}d\sigma\\
&=\int_{\Sigma}[|\nabla\phi|^2-(|A|_g^{2}+\overline{\textrm{Ric}}_{f}(\nu,\nu))\phi^2]e^{-f}d\sigma.
\end{split}
\end{equation}
The $L_{f}$-index, denoted by $L_{f}$-ind, of $(\Sigma,  \bar{g})$  is defined to be the maximum of the dimensions of negative definite  subspaces of $B_f$.
Hence $(\Sigma,  \bar{g})$ is $L_{f}$-stable if and only if $L_{f}$-ind$=0$. Clearly the definition of $L_f$-index is equivalent to the corresponding definition using the variational field $T$ before.

Also, for minimal hypersurface $i: (\Sigma, \tilde{g})\to (M^{n+1}, \tilde{g})$,  it is well known that
if $\psi$ is a compactly supported smooth function on $\Sigma$, then the second variational  of the volume $\tilde{V}$ of $(\Sigma, i^*\tilde{g})$ is given by
\begin{equation}
\frac{d^{2}}{dt^{2}}\tilde{V}(\Sigma_{t})\biggr|_{t=0}=-\int_{\Sigma}\psi J(\psi)d\tilde{\sigma},
\end{equation}
where  $\tilde{A}$ denotes the second fundamental form of $(\Sigma, \tilde{g})$,  $\tilde{\nu}$ denotes the unit normal of $(\Sigma, \tilde{g})$, and $J=\triangle_{\tilde{g}}+|\tilde{A}|^{2}_{\tilde{g}}+\widetilde{\text{Ric}}(\tilde{\nu},\tilde{\nu})$ is stability operator (or  the Jacobi operator) of $(\Sigma, \tilde{g})$.

It  holds that, from (\ref{JL}),
\begin{prop} Let $(\Sigma^{n}, g)$ be an $f$-minimal hypersurface immersed in $(M, \bar{g})$. Then  for all $\phi\in\mathcal{C}^{\infty}_{o}(\Sigma)$,
\begin{equation}\label{JL1}
\int_{\Sigma}(e^{-\frac{f}{n}}\phi)J(e^{-\frac{f}{n}}\phi)e^{-f}d\sigma=\int_{\Sigma}\phi L_f(\phi)e^{-f}d\sigma.
\end{equation}
\end{prop}

\begin{cor}\label{JLf} For $\phi\in C^{\infty}(\Sigma)$,  $J(e^{-\frac{f}{n}}\phi)=e^{\frac{f}{n}}L_f(\phi).$
\end{cor}
\begin{cor}  $L_f$-ind of $(\Sigma,\bar{g})$ is equal to the index of $(\Sigma, \tilde{g})$. In particular,
 $(\Sigma, \bar{g})$ is $L_f$-stable if and only if $(\Sigma, \tilde{g})$ is stable in   $(M, \tilde{g})$.
\end{cor}


\end{document}